\documentclass[12pt]{article}

\usepackage{color}

\input{preamble.tex}

\begin{document}

\title{A note on the role of projectivity in likelihood-based inference for random graph models}

\author{M.\ Schweinberger\footnote{The first two authors made equal contributions.}
\and P.\ N.\ Krivitsky
\and C.\ T.\ Butts}

\date{}

\maketitle

\begin{abstract}
There is widespread confusion about the role of projectivity in likelihood-based inference for random graph models.
The confusion is rooted in claims that projectivity,
a form of marginalizability,
may be necessary for likelihood-based inference and consistency of maximum likelihood estimators.
We show that likelihood-based superpopulation inference is not affected by lack of projectivity and that projectivity is not a necessary condition for consistency of maximum likelihood estimators. 
\end{abstract}

\section{Introduction}
\label{intro}

\subsection{Motivation}

In the past decade,
network data have attracted much attention and so have models of network data,
known as random graph models \citep[][]{Fien12,HuKrSc12}.
Despite recent advances,
there is widespread confusion about key issues of statistical inference for random graph models.
Chief among them is the role of projectivity, 
a form of marginalizability,
in likelihood-based inference for random graph models.
Based on the work of \citet{ShRi11},
many statisticians have expressed concern that likelihood-based inference for non-projective random graph models may be problematic and their maximum likelihood estimators may be inconsistent.
Since many random graph models are non-projective,
including sparse Bernoulli random graph models \citep{ErRe60} and other random graph models in common use \citep{HuKrSc12},
it is important to clarify the role of projectivity in likelihood-based inference for random graph models. 

\subsection{Goal} 

We clarify that likelihood-based superpopulation inference is not affected by lack of projectivity and that projectivity is not a necessary condition for consistency of maximum likelihood estimators.
In addition,
we argue that consistency under replication may be preferable to consistency under network growth and that consistency under replication does not require projectivity of random graph models.

\subsection{Projectivity}
\label{intro1}

Consider an exponential family of distributions $\{\M_{\verts,\cnmap},\, \cnmap\in\Xi\}$ for a random graph with a set of nodes $\verts$ and a set of edges $\mathscr{E} \subset \verts \times \verts$,
where $\cnmap \in \Xi$ is the natural parameter vector of the exponential family.
The natural parameter vector $\cnmap \equiv \cnmap(\theta, \verts)$ may be a function of a parameter vector $\theta \in \Theta$ and $\verts$.
An example are the classic Bernoulli$(\denspar)$ random graph models \citep{ErRe60},
which assume that edges are independent and identically distributed Bernoulli$(\denspar)$ random variables, 
with $\denspar$ denoting the probability of an edge.
Bernoulli$(\denspar)$ random graph models are exponential-family random graph models with the number of edges as sufficient statistic and natural parameter $\eta = \logit(\denspar)$.
Some of the most interesting random graph phenomena occur in the sparse graph regime where $\denspar_{|\verts|}$ and hence $\eta(\denspar_{|\verts|}) = \logit(\denspar_{|\verts|})$ depend on the size $|\verts|$ of $\verts$ \citep{ErRe60}.
For example,
the threshold for connectivity of Bernoulli$(\denspar_{|\verts|})$ random graphs corresponds to $\denspar_{|\verts|} = (\log |\verts|)\, /\, |\verts|$ \citep{ErRe60},
which implies that the natural parameter $\eta(\denspar_{|\verts|}) = \logit(\denspar_{|\verts|})$ depends on $|\verts|$.

In a widely read paper,
\citet{ShRi11} defined projectivity of exponential-family random graph models as follows. 
Let $\verts'\subset\verts$ be a subset of nodes and $\M_{\verts\to\verts',\cnmap(\theta, \verts)}$ be the distribution of the subgraph induced by $\verts'$,
that is,
the marginalization of $\M_{\verts,\cnmap(\theta, \verts)}$ with respect to edge variables involving nodes in $\verts\setminus\verts'$.
An exponential-family random graph model is projective if $\M_{\verts',\cnmap(\theta, \verts')} = \M_{\verts\to\verts',\cnmap(\theta, \verts)}$ and $\cnmap(\theta, \verts') = \cnmap(\theta, \verts)$ for all $\theta \in \Theta$  and all $\verts' \subset \verts$,
regardless of the size of $\verts'$.
For example, 
Bernoulli$(\denspar_{|\verts|})$ random graph models are projective as long as $\denspar_{|\verts|}$ and $\eta(\denspar_{|\verts|}) = \logit(\denspar_{|\verts|})$ do not depend on $\abs{\verts}$.
When $\denspar_{|\verts|}$ and $\eta(\denspar_{|\verts|}) = \logit(\denspar_{|\verts|})$ do depend on $\abs{\verts}$, 
Bernoulli$(\denspar_{|\verts|})$ random graph models are not projective.
It follows that Bernoulli$(\denspar_{|\verts|})$ random graph models are not projective in the sparse graph regime,
where some of the most interesting random graph phenonema occur \citep{ErRe60}.
Many other random graph models are likewise non-projective,
in part because random graphs may be sparse \citep{ErRe60} and in part because edges may be dependent random variables \citep{ShRi11}.

\subsection{Projectivity and statistical inference}
\label{statistical.inference}

\citet[][page~509]{ShRi11} assumed that researchers \emph{``fit ERGMs (by maximum likelihood or pseudo-likelihood) to the observed sub-network, and then extrapolate the same model, with the same parameters, to the whole network''},
where ``ERGMs'' refers to exponential-family random graph models. 
\citet[][page~510]{ShRi11} went on to argue that,
when random graph models are non-projective,
\emph{``the parameter estimates obtained from a sub-network may not provide reliable estimates of... the parameters of the whole network, rendering the task of statistical inference based on a sub-network ill-posed.''}  
Here,
\citeauthor{ShRi11} make the important point that,
given an observed subgraph $\ooo$ of a graph $\net_{\verts}$ with $\verts'\subset\verts$, 
naive statistical inference for $\M_{\verts,\cnmap(\theta, \verts)}$ based on $\M_{\verts',\cnmap(\theta, \verts')}$ may be problematic.

While interesting,
the results of \citeauthor{ShRi11} have been widely misinterpreted as implying that consistent estimation of non-projective random graph models may not be possible.
For example,
\citet[page~831]{Fien12} writes:
\emph{\enquote{The Shalizi-Rinaldo results also explain the sense in which one can or cannot get the consistency of maximum likelihood estimation for ERGMs.}}
\citet{Fien12} seems to suggest that projectivity is necessary for consistency of maximum likelihood estimators and that consistent estimation of non-projective random graph models may hence not be possible.
Others have voiced similar concerns,
sometimes in writing \citep[][]{YaLeZh11} and more often in personal communications and professional meetings.

\section{The likelihood is not affected by lack of projectivity}
\label{super}

The motivating example of \citet{ShRi11} concerns likelihood-based superpopulation inference.
In other words,
the goal is to infer the population model $\M_{\verts,\cnmap(\theta, \verts)}$ that generated a population graph $\net_{\verts}$ defined on a finite population of nodes $\verts$,
where $\net_{\verts}$ is unobserved but a subgraph $\net_{\verts'}$ of $\net_{\verts}$ induced by a subset of nodes $\verts' \subset \verts$ is observed.

One problem, 
which has been the source of considerable confusion,
is that\linebreak 
\citeauthor{ShRi11} considered statistical inference based on $\M_{\verts',\cnmap(\theta, \verts')}(\OOO = \ooo)$,
despite the fact that the likelihood is not proportional to $\M_{\verts',\cnmap(\theta, \verts')}(\OOO = \ooo)$ unless the population model and sampling design satisfy additional conditions. 
In general,
the likelihood is proportional to the probability of the observed data \citep{Fi22}.
In particular,
if a population graph $\net_{\verts}$ is generated by population model $\M_{\verts,\cnmap(\theta, \verts)}$ and a subgraph $\net_{\verts'}$ of $\net_{\verts}$ induced by a subset of nodes $\verts' \subset \verts$ is sampled by an ignorable sampling design \citep{HaGi09},
then the likelihood is
\begin{equation}
\label{likelihood}
\lik(\theta;\, \net_{\verts'})
\;\propto\; {\sum_{\net_{\verts}\, \in\, \netspace_{\verts}(\ooo)}}\M_{\verts,\cnmap(\theta, \verts)}(\Net_{\verts} = \net_{\verts}),
\end{equation}
where $\netspace_{\verts}(\ooo)$ is the set of all graphs on $\verts$ whose induced subgraph on $\verts'$ is $\ooo$.

Two conclusions follow.
First,
by construction,
the likelihood \eqref{likelihood} is proportional to the marginalization ${\sum_{\net_{\verts}\, \in\, \netspace_{\verts}(\ooo)}}\M_{\verts,\cnmap(\theta, \verts)}(\Net_{\verts} = \net_{\verts})$ and is hence not affected by lack of projectivity of $\M_{\verts,\cnmap(\theta, \verts)}$. 
Second,
the misspecified likelihood $\M_{\verts',\cnmap(\theta, \verts')}(\OOO = \ooo)$ of \citeauthor{ShRi11} is not, 
in general, 
proportional to the proper likelihood \eqref{likelihood}, 
hence the results of \citeauthor{ShRi11} are not pertinent to likelihood-based superpopulation inference.

\section{Projectivity is not necessary for consistency}
\label{sec:meandeg}

\citet{Fien12} and others suggested that projectivity may be necessary for consistency of maximum likelihood estimators.
We demonstrate that projectivity is not necessary for consistency of maximum likelihood estimators by a counterexample.  

Consider a sequence of classic Bernoulli$(\denspar_{|\verts|})$ random graphs \citep{ErRe60} with $|\verts|$ nodes and size-dependent edge probabilities $\denspar_{|\verts|}$,
where $|\verts| = 1, 2, \dotsc$.\linebreak
\citet{KrHaMo11} proposed the parameterization $\denspar_{|\verts|} = \ilogit(\param-\log |\verts|)$.
Here,
the probability of an edge $\denspar_{|\verts|}$ depends on a size-invariant parameter $\param \in \mathbb{R}$ and a size-dependent offset $\log {|\verts|}$,
and so does the natural parameter $\eta(\theta, \verts) = \logit(\denspar_{|\verts|}) = \theta - \log |\verts|$.
This parameterization is motivated by invariance considerations:
the expected number of edges of each node tends to $\exp(\param)$ as $|\verts|\to\infty$ and is hence invariant to network size
\citep{KrHaMo11}.
Such models are not projective,
because $\denspar_{|\verts|}$ and $\eta(\theta, \verts) = \theta - \log |\verts|$ depend on $|\verts|$.
Despite the lack of projectivity,
the maximum likelihood estimator $\widehat\param_{|\verts|}$ of the size-invariant parameter $\theta$ is a consistent estimator of $\theta$ as ${|\verts|}\to\infty$ \citep[][Theorem 3.1]{KrKo14}.  
Therefore,
projectivity is not necessary for consistency of maximum likelihood estimators.

\section{Consistency under replication does not require projectivity}

Consistency under network growth,
as considered by \citet{ShRi11} and others, 
may not be desirable in the first place, 
because the size of many networks is bounded above and networks of different sizes are governed by different substantive processes. 
Consistency under replication of similar-sized graphs from a common generating process may be preferable to consistency under network growth.
For example,
consistency results may be obtained when $N$ independent graphs of the same size $\net_{\verts}^{(1)}, \ldots,\net_{\verts}^{(N)}$ from $\M_{\verts,\cnmap(\theta, \verts)}$ are observed and $N \to \infty$.
Consistency under replication does not require projectivity of $\M_{\verts,\cnmap(\theta, \verts)}$.

If it is not possible to observe independent graphs of the same size,
consistency under replication is possible when a graph consists of subgraphs of similar size from a common generating process \citep{ScHa13,ScSt16}.
An example is a friendship network of high school students,
where the subgraphs of similar size correspond to friendship networks within and between high schools of similar size.

\section{Conclusion}

Many real-world network processes are not believed to be projective,
because networks of different sizes are governed by different substantive processes.
Thus,
superimposing projectivity on random graph models may be undesirable.
Indeed,
projectivity is not necessary for likelihood-based inference.

\hide{

\section*{Acknowledgements} 

We are grateful to Noel A.\ C.\ Cressie, Mathias Drton, Mark S.\ Handcock, Eric D.\ Kolaczyk, Martina Morris, and Thomas Richardson for helpful discussions.
This research was supported in part by the National Science Foundation and the Army Research Office,
U.S.A.

}



\end{document}